# Gamma, Psi, Bernoulli Functions via Hurwitz Zeta Function.

VIVEK V. RANE

A-3/203, Anand Nagar,

Dahisar (East)

Mumbai – 400 068

INDIA

e-mail address : v_v_rane@ yahoo.co.in

**Abstract** : Using three basic facts concerning Hurwitz zeta function, we give new natural proofs of the known results on Bernoulli polynomials, gamma function and also obtain Gauss' expression for psi function at a rational point all in a unified fashion. We also give a new proof of the relation between log gamma and the derivative of the Hurwitz zeta function including that of Stirling's expression for

 log gamma .

# Gamma , Psi , Bernoulli Functions via Hurwitz Zeta Function

VIVEK V . RANE

A-3/203 , Anand Nagar ,

Dahisar (East)

Mumbai – 400 068

INDIA

e-mail address : v_v_rane@ yahoo.co.in

The object of this paper is to indicate that the theory of Euler's gamma function , that of Gauss ' psi function and that of Bernoulli functions can be developed naturally , with ease from the theory of Hurwitz's zeta function and its derivatives .

For a complex variable s and for a complex number α≠0 , -1, -2 , ………. let ζ(s,α ) be the Hurwitz's zeta function defined by $\zeta(s,\alpha)=\sum_{n=0}^{\infty}(n+\alpha)^{-s}$ for Re s>1 ; and its analytic continuation . Let ζ(s,1)=ζ(s) be the Riemann zeta function . If χ(mod q) is a Dirichlet character modulo an integer q≥1 and if L(s,χ) is the corresponding Dirichlet L-series , then in view of L(s,χ)=$q^{-s}\sum_{a=1}^{q}\chi(a)\zeta(s,\frac{a}{q})$ , it is already clear from author [2] and many other works of the author that the theory of Dirichlet L-series can be developed naturally from the theory of ζ(s,α) .

For integers n≥0 , let $B_n(\alpha)$ be the Bernoulli polynomials in α defined by

$\frac{ze^{\alpha z}}{e^z-1} = \sum_{n=0}^{\infty} B_n(\alpha)\frac{z^n}{n!}$ for |z|<2π ; and let $B_n=B_n(0)$ be Bernoulli numbers . Note that $B_0(\alpha)$=1 . Then it is known that if n≥0 is an integer , $B_{n+1}(\alpha)$=(-n-1)ζ(-n,α) and this fact has been proved in author [4] by comparing the Fourier series expressions of ζ(-n,α) and $B_{n+1}(\alpha)$ . As a function of s , ζ(s,α) is a meromorphic function with a simple pole at s=1 with residue 1 . As a function of α , ζ(s,α) is an analytic function of α except possibly for α=0,-1,-2 , ………. However , if m≥0 is an integer , then ζ(-m,α) is a polynomial in α and thus is an entire function of α (See author [2]) . For any integer r≥0 , we write

$\zeta^{(r)}(s,\alpha) = \frac{\partial^r}{\partial s^r} \zeta(s,\alpha)$. In particular, $\zeta'(s,\alpha) = \frac{\partial}{\partial s}\zeta(s,\alpha)$. As a function of α, $\zeta^{(r)}(s,\alpha)$ is an analytic function of α except possibly for α=0, -1, -2, ………. (See author [2], author [3].) If s≠1 and α≠0, -1, -2, ………, then $\frac{\partial^{r_1}}{\partial \alpha^{r_1}} \frac{\partial^{r_2}}{\partial s^{r_2}} \zeta(s,\alpha) = \frac{\partial^{r_2}}{\partial s^{r_2}} \frac{\partial^{r_1}}{\partial \alpha^{r_1}} \zeta(s,\alpha)$, where $r_1, r_2 \geq 0$ are integers. See for example author [3].

In particular, if s≠1 and if α≠0, -1, -2, ………, then $\frac{\partial}{\partial \alpha} \frac{\partial}{\partial s} \zeta(s,\alpha) = \frac{\partial}{\partial s} \frac{\partial}{\partial \alpha} \zeta(s,\alpha)$.

Next, let Γ(α) be the Euler's gamma function defined by Weierstrass' product namely

$\frac{1}{\Gamma(\alpha)} = \alpha e^{\gamma \alpha} \prod_{n \geq 1}(1+\frac{\alpha}{n})e^{-\frac{\alpha}{n}}$, where γ is Euler's constant. Then it is known that ζ'(0,α)=log$\frac{\Gamma(\alpha)}{\sqrt{2\pi}}$. In what follows, this fact will be proved with ease afresh as Proposition 5, wherein we shall also obtain Stirling's expression for log Γ(α). Next, let ψ(α)=$\frac{d}{d\alpha}$ log Γ(α)=$\frac{d}{d\alpha}$ log $\frac{\Gamma(\alpha)}{\sqrt{2\pi}}$ =$\frac{\partial}{\partial \alpha}$ζ'(0,α). In Proposition 4 below, we shall obtain Gauss' expression for $\psi\left(\frac{a}{q}\right)$, where 1≤a<q are integers.

We shall show that by defining

1) $B_{n+1}(\alpha) = (-n-1)\zeta(-n,\alpha)$ for n≥ -1

2) log Γ(α)=$\frac{\partial}{\partial s}$ζ(s,α) at s=0

3) ψ(α)=$\frac{\partial}{\partial \alpha} \frac{\partial}{\partial s} \zeta(s,\alpha)$ at s=0

we can develop naturally the theory of Bernoulli polynomials, gamma function and Gauss' ψ function in terms of ζ(s,α) and its derivatives as in the case of Dirichlet L-series.

Basically, there are three facts concerning ζ(s,α).

1) $\frac{\partial}{\partial \alpha} \zeta(s,\alpha) = -s\zeta(s+1,\alpha)$

2) $\zeta(s,\alpha) - \zeta(s,\alpha+1) = \alpha^{-s}$

3) For Re s<1 and for 0<α<1 , ζ(s,α)=$2^s \pi^{s-1} \Gamma(1-s) \sum_{n\geq 1} \sin\left(\frac{\pi s}{2} + 2\pi n\alpha\right) n^{s-1}$

=$\Gamma(1-s)\sum_{|n|\geq 1} e^{2\pi i n\alpha} (2\pi i n)^{s-1}$

Fact I) resulted in the Taylor series expansion of ζ(s,α+1) in author [1] in the form

$\zeta(s,\alpha+1) = \zeta(s) + \sum_{n\geq 1} \frac{(-\alpha)^n}{n!} s(s+1) \ldots (s+n-1)\zeta(s+n)$ and consequently in author [2] , we find that for integral m≥0 , ζ(-m,α)=$\sum_{k=0}^{m} \binom{m}{k} \zeta(-k)\alpha^{m-k} + \alpha^m - \frac{\alpha^{m+1}}{m+1}$ . In author [4] , it has been shown that the above expression for ζ(-m,α) with m≥0 , is equivalent to $B_{m'}(\alpha)=\sum_{k=0}^{m'} \binom{m'}{k} B_k \alpha^{m'-k}$ with $m' = m+1$ . Thus fact II) and fact III) remain to be exploited . Using these facts , we prove with ease the following four propositions .

**Proposition 1** We have

1) $B_m(\alpha+1) - B_m(\alpha) = m\alpha^{m-1}$ for m≥0 .

2) $\Gamma(\alpha+1) = \alpha\Gamma(\alpha)$

**Corollary** : ψ(α+1)=ψ (α) + $\frac{1}{\alpha}$

**Proposition 2** : We have

1) $B_m(1-\alpha) = (-1)^m B_m(\alpha)$

2) $\Gamma(\alpha)\Gamma(1-\alpha) = \frac{\pi}{\sin\pi\alpha}$

**Corollary** : ψ(1-α)=ψ (α) +π cot πα

**Proposition 3** : We have for any integer m≥2 and for any integer n≥1 ,

1) $B_n(\alpha) + B_n\left(\alpha + \frac{1}{m}\right) + B_n(\alpha + \frac{2}{m}) + \cdots \ldots \ldots + B_n\left(\alpha + \frac{m-1}{m}\right) = m^{1-n} B_n(m\alpha)$

2) $\Gamma(\alpha)\Gamma\left(\alpha + \frac{1}{n}\right) \ldots \ldots \ldots \Gamma\left(\alpha + \frac{n-1}{n}\right) = (2\pi)^{\frac{n-1}{2}} n^{\frac{1}{2}-n\alpha} \; \Gamma(n\alpha)$

**Corollary** : We have $\prod_{k=0}^{m-1} \sin\pi(\alpha + \frac{k}{m}) = 2^{1-m} \sin\pi m\alpha$.

Note that in the above three propositions, the results are in pairs. The results pertaining to Bernoulli polynomials are about sums and the results pertaining to gamma function are about products.

Next, we state remaining Propositions.

**Proposition 4** : We have If 1≤a<q are integers,

then $\psi\left(\frac{a}{q}\right) = -(\gamma + \log q) + \sum_{r=1}^{q-1} \cos\frac{2\pi r a}{q} \log 2\sin\frac{\pi r}{q} + \frac{\pi}{q}\sum_{r=1}^{q-1} r \sin\frac{2\pi r a}{q}$

Next with Weierstrass product definition of $\Gamma(\alpha)^{-1}$ and using Euler's summation formula, we state and prove Proposition 5. In what follows, we shall write $\phi(u) = u - [u] - \frac{1}{2}$, where [u] denotes integral part of the real variable u.

**Proposition 5** : Let $\phi(u) = u - [u] - \frac{1}{2}$.

Then log Γ(α)=(α-1/2) log α − α +1 + $\int_1^\infty \frac{\phi(u)}{u} du - \int_0^\infty \frac{\phi(u)}{u+\alpha} du$ and $\zeta'(0, \alpha) = \log\frac{\Gamma(\alpha)}{\sqrt{2\pi}}$.

Next, we prove our propositions. We shall prove Proposition 5 first.

**Proof of Proposition 5** : From Weierstrass product expression for $\frac{1}{\Gamma(\alpha)}$ on taking logarithm,

we have $\log \Gamma(\alpha) = \sum_{n\geq 1}(\frac{\alpha}{n} - \log(1 + \frac{\alpha}{n})) - \gamma\alpha - \log\alpha$.

On using Euler's summation formula

$$\log \Gamma(\alpha) = \frac{1}{2}(\alpha - \log(1+\alpha)) - \gamma\alpha - \log \alpha + \int_1^\infty (\frac{\alpha}{u} - \log(1+\frac{\alpha}{u}))du +$$

$$+ \int_1^\infty \emptyset(u) \frac{d}{du}(\frac{\alpha}{u} - \log(1+\frac{\alpha}{u}))du , \text{ where } \phi(u) = u - [u] - 1/2 .$$

Next $\int_1^\infty (\frac{\alpha}{u} - \log(1+\frac{\alpha}{u}))du = \left[(u+\alpha)\log\frac{u}{u+\alpha}\right]_{u=1}^\infty = \lim_{u \to \infty}(u+\alpha)\log\frac{u}{u+\alpha} - (1+\alpha)\log\frac{1}{1+\alpha}$

$= -\alpha + (1+\alpha)\log(1+\alpha)$ .

Thus $\log \Gamma(\alpha) = (\alpha + \frac{1}{2}) \log(1+\alpha) - (\gamma + \frac{1}{2})\alpha - \log \alpha + \int_1^\infty \phi(u)\frac{d}{du}\left(\frac{\alpha}{u} - \log\left(1+\frac{\alpha}{u}\right)\right)du$ .

Next , $\int_1^\infty \phi(u)\frac{d}{du}\left(\frac{\alpha}{u} - \log(1+\frac{\alpha}{u})\right)du = \alpha \int_1^\infty \frac{\emptyset(u)}{u}\left(\frac{1}{u+\alpha} - \frac{1}{u}\right)du$

$= \int_1^\infty \phi(u)\left(\frac{1}{u} - \frac{1}{u+\alpha}\right)du - \alpha \int_1^\infty \frac{\emptyset(u)}{u^2}du = \int_1^\infty \frac{\emptyset(u)}{u}du - \int_1^\infty \frac{\emptyset(u)}{u+\alpha}du + \alpha\left(\gamma - \frac{1}{2}\right)$

$= \int_1^\infty \frac{\emptyset(u)}{u}du - \int_0^\infty \frac{\emptyset(u)}{u+\alpha}du + \alpha\left(\gamma - \frac{1}{2}\right) + 1 + (\alpha+\frac{1}{2})\log\alpha - (\alpha+\frac{1}{2})\log(1+\alpha)$ .

Thus $\log \Gamma(\alpha) = (\alpha - \frac{1}{2})\log\alpha - \alpha + 1 + \int_1^\infty \frac{\emptyset(u)}{u}du - \int_0^\infty \frac{\emptyset(u)}{u+\alpha}du$ .

This gives an expression for $\log \Gamma(\alpha)$ . Note that integrating by parts repeatedly , $\int_0^\infty \frac{\emptyset(u)}{u+\alpha}$ can be developed into an asymptotic series in $\alpha$ .

Next , we shall prove $\zeta'(0, \alpha) = \log\frac{\Gamma(\alpha)}{\sqrt{2\pi}}$ . For $0<\alpha\leq 1$ , Re $s>1$ and for arbitrary $x>0$ ,

we have $\zeta(s, \alpha) = \sum_{n\geq 0}(n+\alpha)^{-s} = \sum_{0\leq n \leq x-\alpha}(n+\alpha)^{-s} + \sum_{n>x-\alpha}(n+\alpha)^{-s}$ .

We use Euler's summation for $\sum_{n>x-\alpha}(n+\alpha)^{-s}$ .

This gives for x>0 , 0<α≤1 and for Re s>0 ,

$$\zeta(s,\alpha) = \sum_{0 \leq n \leq x-\alpha}(n+\alpha)^{-s} + \frac{x^{1-s}}{s-1} + (x-\alpha-[x-\alpha]-\frac{1}{2})x^{-s} - s\int_x^\infty \frac{\phi(u-\alpha)}{u^{s+1}} du .$$

Choosing x=α , we get $\zeta(s,\alpha) = \alpha^{-s} + \frac{\alpha^{1-s}}{s-1} - \frac{\alpha^{-s}}{2} - s\int_\alpha^\infty \frac{\phi(u-\alpha)}{u^{s+1}} du$

so that $\zeta(s,\alpha) = \frac{\alpha^{-s}}{2} + \frac{\alpha^{1-s}}{s-1} - s\int_\alpha^\infty \frac{\phi(u-\alpha)}{u^{s+1}} du$ .

Thus $\lim_{s \to 0} \frac{\zeta(s,\alpha) - \frac{\alpha^{-s}}{2} - \frac{\alpha^{1-s}}{s-1}}{s} = - \lim_{s \to 0} \int_\alpha^\infty \frac{\phi(u-\alpha)}{u^{s+1}} du$ .

As s→ 0 , $\zeta(s,\alpha) - \frac{\alpha^{-s}}{2} - \frac{\alpha^{1-s}}{s-1} = \zeta(0,\alpha) - \frac{1}{2} + \alpha = (\frac{1}{2}-\alpha) - \frac{1}{2} + \alpha \to 0$ .

Thus the left hand side is of the form $\frac{0}{0}$ as s→ 0 through real values . Hence by L'Hospital's rule ,

We have $\lim_{s \to 0} \frac{\zeta(s,\alpha) - \frac{\alpha^{-s}}{2} - \frac{\alpha^{1-s}}{s-1}}{s} = \lim_{s \to 0} \left(\zeta'(s,\alpha) + \frac{\alpha^{-s}}{2} \log\alpha + \frac{\alpha^{1-s} \log\alpha}{s-1} + \frac{\alpha^{1-s}}{(s-1)^2}\right)$

$= \zeta'(0,\alpha) + \frac{1}{2}\log\alpha - \alpha\log\alpha + \alpha$ and $\lim_{s \to 0} \int_\alpha^\infty \frac{\phi(u-\alpha)}{u^{s+1}} du = \int_\alpha^\infty \frac{\phi(u-\alpha)}{u} du$ .

Thus we have $\zeta'(0,\alpha) + \frac{1}{2}\log\alpha - \alpha\log\alpha + \alpha = -\int_\alpha^\infty \frac{\phi(u-\alpha)}{u} du$ .

That is $\zeta'(0,\alpha) + \frac{1}{2}\log\alpha - \alpha\log\alpha + \alpha = -\int_0^\infty \frac{\phi(u)}{u+\alpha} du$ .

That is $\zeta'(0,\alpha) = (\alpha - \frac{1}{2})\log\alpha - \alpha - \int_0^\infty \frac{\phi(u)}{u+\alpha} du$ .

This gives $\zeta'(0,\alpha) = \log \Gamma(\alpha) + c$  with c  a constant , on comparing with the above expression for

$\log \Gamma(\alpha)$ . Putting α=1 , we get c=ζ'(0)= $-\frac{1}{2}$ log 2π .

Thus we have $\zeta'(0, \alpha) = \log \frac{\Gamma(\alpha)}{\sqrt{2\pi}}$.

**Proof of Proposition 1**: I) We have $\zeta(s, \alpha) - \zeta(s, \alpha + 1) = \alpha^{-s}$.

Taking s= - m , where m≥0 is an integer , we have $-\frac{B_{m+1}(\alpha)}{m+1} + \frac{B_{m+1}(\alpha+1)}{m+1} = \alpha^m$.

That is , $B_{m+1}(\alpha + 1) - B_{m+1}(\alpha) = (m + 1)\alpha^m$

II) From $\zeta(s, \alpha) - \zeta(s, \alpha + 1) = \alpha^{-s}$ ,

on differentiation with respect to s , we have $\zeta'(s, \alpha) - \zeta'(s, \alpha + 1) = -\alpha^{-s} \log \alpha$.

Putting s=0 , we have $\zeta'(0, \alpha) - \zeta'(0, \alpha + 1) = - \log \alpha$.

Thus $\log \frac{\Gamma(\alpha)}{\sqrt{2\pi}} - \log \frac{\Gamma(\alpha+1)}{\sqrt{2\pi}} = $ - log α , in view of Proposition 5 .

That is log Γ(α+1) – log Γ(α) = log α  or  Γ(α+1) = αΓ(α) .

**Proof of Corollary** : Differentiating the expression  log Γ(α+1) – log Γ(α) = log α ,

we get ψ(α+1)-ψ(α) = $\frac{1}{\alpha}$ .

**Proof of Proposition 2** : I) We have for  Re s<1 and for 0<α<1 ,

$\zeta(s, \alpha) = 2^s \pi^{s-1} \Gamma(1-s) \sum_{n\geq 1} \sin\left(\frac{\pi s}{2} + 2\pi n \alpha\right) n^{s-1}$

So that $\zeta(s, 1 - \alpha) = 2^s \pi^{s-1} \Gamma(1-s) \sum_{n\geq 1} \sin\left(\frac{\pi s}{2} - 2\pi n \alpha\right) n^{s-1}$.

Writing s= -m , where m≥0 is an integer ,

we have $\zeta(-m, 1 - \alpha) = -2^{-m} \pi^{-m-1} \Gamma(1+m) \sum_{n\geq 1} \sin\left(\frac{\pi m}{2} + 2\pi n \alpha\right) n^{-m-1}$

whereas $\zeta(-m, \alpha) = 2^{-m}\pi^{-m-1}\Gamma(1+m) \sum_{n\geq 1} \sin\left(\frac{-\pi m}{2} + 2\pi n\alpha\right) n^{-m-1}$

This gives $\zeta(-m, 1-\alpha) = (-1)^{m+1} \zeta(-m, \alpha)$

That is, $B_{m+1}(1-\alpha) = (-1)^{m+1} B_{m+1}(\alpha)$.

II) We have for Re s<1 and for 0<α<1, $\zeta(s, \alpha) = 2(2\pi)^{s-1}\Gamma(1-s) \sum_{n\geq 1} \sin\left(\frac{\pi s}{2} + 2\pi n\alpha\right) n^{s-1}$

Consider $\zeta(s, 1-\alpha) = 2(2\pi)^{s-1}\Gamma(1-s) \sum_{n\geq 1} \sin\left(\frac{\pi s}{2} - 2\pi n\alpha\right) n^{s-1}$.

Differentiating with respect to s,

$\zeta'(s, 1-\alpha) = 2(2\pi)^{s-1}\{(\log 2\pi \cdot \Gamma(1-s) - \Gamma'(1-s)) \sum_{n\geq 1} \sin\left(\frac{\pi s}{2} - 2\pi n\alpha\right) \cdot n^{s-1}$

$+\Gamma(1-s) \sum_{n\geq 1} \left(\frac{\pi}{2}\cos\left(\frac{\pi s}{2} - 2\pi n\alpha\right) + \sin\left(\frac{\pi s}{2} - 2\pi n\alpha\right) \cdot \log n\right) n^{s-1}\}$

Thus $\zeta'(0, 1-\alpha) = \frac{1}{\pi}\{(\log 2\pi + \gamma) \sum_{n\geq 1} \frac{(-\sin 2\pi n\alpha)}{n} + \sum_{n\geq 1} \frac{(\frac{\pi}{2}\cos 2\pi n\alpha - \log n \cdot \sin 2\pi n\alpha)}{n}\}$

Replacing α by 1-α, $\zeta'(0, \alpha) = \frac{1}{\pi}\{(\log 2\pi + \gamma) \sum_{n\geq 1} \frac{\sin 2\pi n\alpha}{n} +$

$+ \sum_{n\geq 1} \frac{1}{n}\left(\frac{\pi}{2}\cos 2\pi n\alpha + \log n \cdot \sin 2\pi n\alpha\right)\}$.

Thus $\zeta'(0,\alpha) + \zeta'(0, 1-\alpha) = \sum_n \frac{\cos 2\pi n\alpha}{n} = \frac{1}{2}\sum_{n\geq 1}\frac{(e^{2\pi i n\alpha} + e^{-2\pi i n\alpha})}{n}$

$= -\frac{1}{2}\left(\log(1 - e^{2\pi i\alpha}) + \log(1 - e^{-2\pi i\alpha})\right) = -\frac{1}{2}\{\log e^{\pi i\alpha}(e^{-\pi i\alpha} - e^{\pi i\alpha}) + \log e^{-\pi i\alpha}(e^{\pi i\alpha} - e^{-\pi i\alpha})\}$

$= -\frac{1}{2}\{\log(e^{\pi i\alpha} - e^{-\pi i\alpha}) + \log(-e^{\pi i\alpha}) + \log(e^{\pi i\alpha} - e^{-\pi i\alpha}) + \log e^{-\pi i\alpha}\}$

$= -\frac{1}{2}\{2\log \frac{(e^{\pi i\alpha} - e^{-\pi i\alpha})}{2i} + \log(-2ie^{\pi i\alpha}) + \log(2ie^{-\pi i\alpha})\}$

$= -\frac{1}{2}\{2 \log \sin \pi\alpha + \log(-4i^2)\} = -(\log \sin \pi\alpha + \log 2) = -\log 2 \sin \pi\alpha$

Thus $\zeta'(0,\alpha) + \zeta'(0, 1-\alpha) + \log 2\sin\pi\alpha = 0$

That is $\log \frac{\Gamma(\alpha)}{\sqrt{2\pi}} + \log \frac{\Gamma(1-\alpha)}{\sqrt{2\pi}} + \log 2\sin\pi\alpha = 0$.

That is, $\log \frac{\Gamma(\alpha)\Gamma(1-\alpha)\sin\pi\alpha}{\pi} = 0$. That is, $\Gamma(\alpha)\Gamma(1-\alpha)\frac{\sin\pi\alpha}{\pi} = 1$

That is $\Gamma(\alpha)\Gamma(1-\alpha) = \frac{\pi}{\sin\pi\alpha}$.

**Proof of Corollary** : From log Γ(α)+log Γ(1-α) = log π-log sin πα ,

on differentiation , we have ψ(α)=ψ(1-α) - π cot πα

**Proof of Proposition 3** : I) For Re s>1 , $\zeta(s,\alpha) = \sum_{n\geq 0}(n+\alpha)^{-s}$.

For an integer m≥2 , consider $\zeta(s, \alpha + \frac{1}{m}) = \sum_{n\geq 0}(n + \alpha + \frac{1}{m})^{-s} = m^s \sum_{n\geq 0}(mn + 1 + m\alpha)^{-s}$

$= m^s\{\sum_{n\geq 0}(n+m\alpha)^{-s} - \sum_{n\equiv 0 (mod\ m)}(n+m\alpha)^{-s} - \sum_{r=2}^{m-1}\sum_{n\equiv r(mod\ m)}(n+m\alpha)^{-s}\}$

$= m^s\{\sum_{n\geq 0}(n+m\alpha)^{-s} - \sum_{k\geq 0}(mk+m\alpha)^{-s} - \sum_{r=2}^{m-1}\sum_{k\geq 0}(mk+r+m\alpha)^{-s}\}$

$= m^s\{\zeta(s, m\alpha) - m^{-s}\sum_{k\geq 0}(k+\alpha)^{-s} - m^{-s}\sum_{r=2}^{m-1}(k+\alpha+\frac{r}{m})^{-s}\}$

$= m^s \zeta(s, m\alpha) - \zeta(s,\alpha) - \sum_{r=2}^{m-1}\zeta(s, \alpha + \frac{r}{m})$

Thus $\sum_{r=0}^{m-1} \zeta(s, \alpha + \frac{r}{m}) = m^s\zeta(s, m\alpha)$

Writing s= - (n-1) , where n≥1 is an integer , we have $\sum_{r=0}^{m-1} B_n(\alpha + \frac{r}{m}) = m^{1-n}B_n(m\alpha)$.

Next differentiating $\sum_{r=0}^{m-1} \zeta(s, \alpha + \frac{r}{m}) = m^s \zeta(s, m\alpha)$ with respect to s,

we have $\sum_{r=0}^{m-1} \zeta'(s, \alpha + \frac{r}{m}) = m^s \zeta'(s, m\alpha) + m^s \log m \cdot \zeta(s, m\alpha)$.

Putting s=0 , we have $\sum_{r=0}^{m-1} \zeta'(0, \alpha + \frac{r}{m}) = \zeta'(0, m\alpha) + \log m \cdot \zeta(0, m\alpha)$.

That is , $\sum_{r=0}^{m-1} \log \frac{\Gamma\left(\alpha + \frac{r}{m}\right)}{\sqrt{2\pi}} = \log \frac{\Gamma(m\alpha)}{\sqrt{2\pi}} + \log m \cdot (\frac{1}{2} - m\alpha)$.

Thus $\sum_{r=0}^{m-1} \log \Gamma(\alpha + \frac{r}{m}) = \log \Gamma(m\alpha) + (m-1) \log \sqrt{2\pi} + \log \sqrt{m} - m\alpha \log m$.

Thus $\log \prod_{r=0}^{m-1} \Gamma(\alpha + \frac{r}{m}) = \log \frac{(2\pi)^{\frac{m-1}{2}} m^{\frac{1}{2}} \Gamma(m\alpha)}{m^{m\alpha}}$.

Thus $\prod_{r=0}^{m-1} \Gamma(\alpha + \frac{r}{m}) = (2\pi)^{\frac{m-1}{2}} m^{\frac{1}{2} - m\alpha} \Gamma(m\alpha)$.

**Proof of Corollary** : Consider for 0<α<1 , $\prod_{k=0}^{n-1} \Gamma(\alpha + \frac{k}{n}) \cdot \prod_{k=0}^{n-1} \Gamma(-\alpha + \frac{k}{n})$.

$= (2\pi)^{n-1} \cdot n^{\frac{1}{2} - n\alpha} \cdot \Gamma(n\alpha) \cdot n^{\frac{1}{2} + n\alpha} \cdot \Gamma(-n\alpha) = (2\pi)^{n-1} \cdot n \, \Gamma(n\alpha) \Gamma(-n\alpha)$

$= -\frac{(2\pi)^{n-1}}{\alpha} \Gamma(n\alpha) \cdot (-n\alpha \cdot \Gamma(-n\alpha)) = -\frac{(2\pi)^{n-1}}{\alpha} \Gamma(n\alpha) \Gamma(1 - n\alpha) = -\frac{(2\pi)^{n-1}}{\alpha} \frac{\pi}{\sin \pi n\alpha}$

On the other hand $\prod_{k=0}^{n-1} \Gamma(\alpha + \frac{k}{n}) \cdot \prod_{k=0}^{n-1} \Gamma(-\alpha + \frac{k}{n}) = \prod_{k=0}^{n-1} \Gamma(\alpha + \frac{k}{n}) \prod_{k=1}^{n} \Gamma(-\alpha + \frac{n-k}{n})$

$= \Gamma(\alpha)\Gamma(-\alpha) \cdot \prod_{k=1}^{n-1} \Gamma(\alpha + \frac{k}{n}) \Gamma(1 - \alpha - \frac{k}{n}) = \Gamma(\alpha)\Gamma(-\alpha) \prod_{k=1}^{n-1} \frac{\pi}{\sin \pi (\alpha + \frac{k}{n})} = -\frac{\pi^{n-1} \Gamma(\alpha) \Gamma(1-\alpha)}{\alpha \prod_{k=1}^{n-1} \sin \pi (\alpha + \frac{k}{n})}$

$= -\frac{\pi^n}{\alpha} \frac{1}{\sin \pi \alpha \cdot \prod_{k=1}^{n-1} \sin \pi (\alpha + \frac{k}{n})} = -\frac{\pi^n}{\alpha} (\prod_{k=0}^{n-1} \sin \pi (\alpha + \frac{k}{n}))^{-1}$.

Thus we have $-\frac{(2\pi)^{n-1}}{\alpha} \cdot \frac{\pi}{\sin \pi n\alpha} = -\frac{\pi^n}{\alpha} \left(\left(\prod_{k=0}^{n-1} \sin \pi (\alpha + \frac{k}{n})\right)\right)^{-1}$

Thus $\frac{2^{n-1}}{\sin \pi n \alpha} = \frac{1}{\left(\prod_{k=0}^{n-1} \sin \pi (\alpha + \frac{k}{n})\right)}$

Thus $\prod_{k=0}^{n-1} \sin \pi (\alpha + \frac{k}{n}) = 2^{1-n} \sin \pi n \alpha$.

**Proof of Proposition 4**: We have $\psi(\alpha) = \frac{d}{d\alpha} \log \Gamma(\alpha) = \left(\frac{\partial}{\partial \alpha} \frac{\partial}{\partial s} \zeta(s,\alpha)\right)_{s=0}$

$= \left(\frac{\partial}{\partial s} \frac{\partial}{\partial \alpha} \zeta(s,\alpha)\right)_{s=0} = \frac{\partial}{\partial s}(-s\zeta(s+1,\alpha))_{s=0} = \left(\frac{\partial}{\partial s}(-(s-1)\zeta(s,\alpha))\right)_{s=1} = \left(\frac{\partial}{\partial s}(1-s)\zeta(s,\alpha)\right)_{s=1}$

Next, we have for Re s<1 and $0<\alpha<1$

$\zeta(s,\alpha) = 2(2\pi)^{s-1} \Gamma(1-s) \sum_{n \geq 1} \sin\left(\frac{\pi s}{2} + 2\pi n \alpha\right) \cdot n^{s-1}$

So that $(1-s)\zeta(s,\alpha) = 2(2\pi)^{s-1} \Gamma(2-s) \sum_{n \geq 1} \sin(\frac{\pi s}{2} + 2\pi n \alpha) n^{s-1}$

This gives for Re s<0, $(1-s)\zeta(s,\frac{a}{q}) = 2(2\pi)^{s-1} \Gamma(2-s) \sum_{n \geq 1} \sin(\frac{\pi s}{2} + 2\pi n \frac{a}{q}) n^{s-1}$

$= 2(2\pi)^{s-1} \Gamma(2-s) \sum_{r=1}^{q} \sin(\frac{\pi s}{2} + 2\pi r \frac{a}{q}) \cdot \sum_{n \equiv r (\bmod\ q)} n^{s-1}$

$= 2(2\pi q)^{s-1} \Gamma(2-s) \sum_{r=1}^{q} \sin(\frac{\pi s}{2} + 2\pi r \frac{a}{q}) \cdot \zeta(1-s, \frac{r}{q})$.

Thus $\frac{\partial}{\partial s}(1-s)\zeta\left(s, \frac{a}{q}\right) = 2(2\pi q)^{s-1}(\log 2\pi q - \Gamma'(2-s)) \cdot \sum_{r=1}^{q} \sin\left(\frac{\pi s}{2} + \frac{2\pi r a}{q}\right) \cdot \zeta\left(1-s, \frac{r}{q}\right)$

$+ 2(2\pi q)^{s-1} \Gamma(2-s) \sum_{r=1}^{q} \{\frac{\pi}{2} \cos(\frac{\pi s}{2} + \frac{2\pi r a}{q}) \cdot \zeta(1-s, \frac{r}{q}) - \sin(\frac{\pi s}{2} + \frac{2\pi r a}{q}) \zeta'(1-s, \frac{r}{q})\}$

$= 2(2\pi q)^{s-1} \left(\log 2\pi q - \frac{\Gamma'}{\Gamma}(2-s)\right) \Gamma(2-s) \cdot \sum_{r=1}^{q} \sin(\frac{\pi s}{2} + \frac{2\pi r a}{q}) \zeta(1-s, \frac{r}{q})$

$+ 2(2\pi q)^{s-1} \Gamma(2-s) \sum_{r=1}^{q} \{\frac{\pi}{2} \cos(\frac{\pi s}{2} + \frac{2\pi r a}{q}) \zeta(1-s, \frac{r}{q}) - \sin(\frac{\pi s}{2} + \frac{2\pi r a}{q}) \zeta'(1-s, \frac{r}{q})\}$.

Next $2(2\pi q)^{s-1}(\log 2\pi q - \frac{\Gamma'}{\Gamma}(2-s)) \cdot \Gamma(2-s) \sum_{r=1}^{q} \sin(\frac{\pi s}{2} + \frac{2\pi r a}{q})\zeta(1-s,\frac{r}{q})$

=$(1-s)(\log 2\pi q - \frac{\Gamma'}{\Gamma}(2-s)) \{2(2\pi q)^{s-1}\Gamma(1-s) \sum_{r=1}^{q} \sin\left(\frac{\pi s}{2} + \frac{2\pi r a}{q}\right)\zeta\left(1-s,\frac{r}{q}\right)\}$

=$(1-s)(\log 2\pi q - \frac{\Gamma'}{\Gamma}(2-s))\zeta(s,\frac{a}{q}) = -(\log 2\pi q - \frac{\Gamma'}{\Gamma}(2-s))((s-1)\zeta(s,\frac{a}{q}))$

$\to -(\log 2\pi q + \gamma) \cdot 1 = -(log 2\pi q + \gamma)$ as $s \to 1$.

As $s \to 1$, consider $2(2\pi q)^{s-1} \Gamma(2-s) \sum_{r=1}^{q}\{\frac{\pi}{2}\cos(\frac{\pi s}{2} + \frac{2\pi r a}{q})\zeta(1-s,\frac{r}{q}) - \sin\left(\frac{\pi s}{2} + \frac{2\pi r a}{q}\right)\zeta'(1-s,\frac{r}{q})\}$.

This is equal to $2\sum_{r=1}^{q}\{\frac{\pi}{2}\cos(\frac{\pi}{2} + \frac{2\pi r a}{q})\zeta(0,\frac{r}{q}) - \sin\left(\frac{\pi}{2} + \frac{2\pi r a}{q}\right)\zeta'(0,\frac{r}{q})\}$

=$2\left(\sum_{r=1}^{q}\left(-\frac{\pi}{2}\right) \sin\frac{2\pi r a}{q} \cdot \left(\frac{1}{2}-\frac{r}{q}\right) - \sum_{r=1}^{q-1} \cos\frac{2\pi r a}{q} \cdot \log\frac{\Gamma\left(\frac{r}{q}\right)}{\sqrt{2\pi}} - \zeta'(o)\right)$

= - $2\{\sum_{r=1}^{q} \frac{\pi}{2} \sin\frac{2\pi r a}{q} \cdot (\frac{1}{2} - \frac{r}{q}) + \sum_{r=1}^{q-1} \cos\frac{2\pi r a}{q} \log\frac{\Gamma\left(\frac{r}{q}\right)}{\sqrt{2\pi}} + \zeta'(o)\}$

= - $2\{-\sum_{r=1}^{q}\frac{\pi}{2}(\sin\frac{2\pi r a}{q}) \cdot \frac{r}{q} + \sum_{r=1}^{q-1} \cos\frac{2\pi r a}{q} \cdot \log\Gamma(\frac{r}{q})\}$

=$\frac{\pi}{q}\sum_{r=1}^{q} r \sin\frac{2\pi r a}{q} - 2\sum_{r=1}^{q-1}\cos\frac{2\pi r a}{q} \cdot \log\Gamma(\frac{r}{q}) = S$, say.

**Case 1**: Let $q \geq 2$ be odd. Then

$S = \frac{\pi}{q}\sum_{r=1}^{q} r \sin\frac{2\pi r a}{q} - 2\sum_{r \leq \frac{q-1}{2}} \cos\frac{2\pi r a}{q}\left(\log\Gamma\left(\frac{r}{q}\right) + \log\Gamma(1-\frac{r}{q})\right)$

=$\frac{\pi}{q}\sum_{r=1}^{q-1} r \sin\frac{2\pi r a}{q} - 2\sum_{r \leq \frac{q-1}{2}} \cos\frac{2\pi r a}{q} \cdot \log\Gamma(\frac{r}{q})\Gamma(1-\frac{r}{q})$

=$\frac{\pi}{q}\sum_{r=1}^{q-1} r \sin\frac{2\pi r a}{q} - 2\sum_{r \leq \frac{q-1}{2}} \cos\frac{2\pi r a}{q} \cdot \log\frac{\pi}{\sin\frac{\pi r}{q}}$

$= \frac{\pi}{q}\sum_{r=1}^{q-1} r \sin\frac{2\pi r a}{q} + 2\sum_{r \leq \frac{q-1}{2}} \cos\frac{2\pi r a}{q} \log \sin\frac{\pi r}{q} - 2 \log \pi \sum_{r \leq \frac{q-1}{2}} \cos\frac{2\pi r a}{q}$

$= \frac{\pi}{q}\sum_{r=1}^{q-1} r \sin\frac{2\pi r a}{q} + \sum_{r \leq q-1} \cos\frac{2\pi r a}{q} \log \sin\frac{\pi r}{q} - \log \pi \cdot \sum_{r \leq q-1} \cos\frac{2\pi r a}{q}$

$= \frac{\pi}{q}\sum_{r=1}^{q-1} r \sin\frac{2\pi r a}{q} + \sum_{r=1}^{q-1} \cos\frac{2\pi r a}{q} \log \sin\frac{\pi r}{q} + \log \pi$

**Case 2** : If q≥2 is even , we have to deal with the additional term corresponding to r=$\frac{q}{2}$ . However , the term corresponding to r=$\frac{q}{2}$ does not contribute to S= $\frac{\pi}{q}\sum_{r=1}^{q-1} r \sin\frac{2\pi r a}{q} + log\pi +$

$+ \sum_{r=1}^{q-1} \cos\frac{2\pi r a}{q} \cdot \log \sin\frac{\pi r}{q}$ and thus S remains unchanged .

Thus $\psi\left(\frac{a}{q}\right) = -(\gamma + \log 2\pi q) + \frac{\pi}{q}\sum_{r=1}^{q-1} r \sin\frac{2\pi r a}{q} + \log \pi + \sum_{r=1}^{q-1} \cos\frac{2\pi r a}{q} \log \sin\frac{\pi r}{q}$

$= -(\gamma + log 2q) + \frac{\pi}{q}\sum_{r=1}^{q-1} r \sin\frac{2\pi r a}{q} + \sum_{r=1}^{q-1} \cos\frac{2\pi r a}{q} \log \sin\frac{\pi r}{q}$

## References

[1] V.V. Rane , Dirichlet L-function and power series for Hurwitz zeta function , Proc. Indian acad.Sci. (Math Sci.) , vol. 103 , No.1, April 1993, pp 27-39 .

[2] V.V. Rane , Instant Evaluation and Demystification of $\zeta(n), L(n,\chi)$ that Euler, Ramanujan Missed-I (arXiv.org website)

[3] V.V. Rane , Instant Evaluation and Demystification of $\zeta(n), L(n,\chi)$ that Euler, Ramanujan Missed-II (arXiv.org website)

[4] V.V. Rane , Instant multiple zeta values at non-positive integers   (arXiv.org website)